\documentclass[letterpaper]{amsart}
\usepackage{amsfonts,amssymb,amscd,amsmath,enumerate,url,verbatim}
\usepackage[all]{xy}
\SelectTips{eu}{}
\xyoption{curve}
\CompileMatrices
\usepackage[colorlinks]{hyperref}
\newcommand{\lra}{\longrightarrow}
\newcommand{\xra}{\xrightarrow}
\newcommand{\thra}{\twoheadrightarrow}
\newcommand{\ov}{\overline}
\newcommand{\col}{\colon}

\newcommand{\fm}{{\mathfrak m}}
\newcommand{\fp}{{\mathfrak p}}
\newcommand{\fq}{{\mathfrak q}}

\newcommand{\depth}{\operatorname{depth}}
\newcommand{\rank}{\operatorname{rank}}
\newcommand{\HH}[2]{\operatorname{H}_{#1}(#2)}
\newcommand{\F}[2]{{\bf F}_{#1}(#2)}
\newcommand{\Tor}[4]{\operatorname{Tor}_{#1}^{#2}(#3,#4)}
\newcommand{\Hom}[3]{\operatorname{Hom}_{#1}(#2,#3)}
\newcommand{\Poi}[3]{\operatorname{P}_{#1}^{#2}(#3)}
\newcommand{\Poig}[4]{\operatorname{P}_{#1}^{#2}(#3,#4)}
\newcommand{\agr}[2][{}]{{{#2}^{\mathsf g}_{#1}}}
\newcommand{\syz}[3][R]{{\Omega^{#1}_{#2}(#3)}}
\newcommand{\Syz}[3]{\Omega^{#1}_{#2}(#3)}

\newcommand{\ds}[1]{\displaystyle{#1}}
\everymath{\displaystyle}

\theoremstyle{plain}
\newtheorem{theorem}{Theorem}[section]
\newtheorem*{Theorem}{Theorem}
\newtheorem{proposition}[theorem]{Proposition}
\newtheorem{lemma}[theorem]{Lemma}
\newtheorem{corollary}[theorem]{Corollary}
\theoremstyle{definition}
\newtheorem*{definition}{Definition}
\newtheorem{chunk}[theorem]{}
\theoremstyle{remark}
\newtheorem{example}[theorem]{Example}
\newenvironment{bfchunk}{\begin{chunk}\textit}{\end{chunk}}

\newtheorem{remark}[theorem]{Remark}
\numberwithin{equation}{theorem}

\newtheorem*{Case1}{Case 1}
\newtheorem*{Case2}{Case 2}
\newtheorem*{Case3}{Case 3}
\newtheorem*{Step1}{Step 1}
\newtheorem*{Step2}{Step 2}
\newtheorem*{Step3}{Step 3}

\begin{document}

\title[Free resolutions]{Free resolutions over short\\ Gorenstein local rings}

\date{20th October, 2009}

\begin{abstract}
 Let $R$ be a local ring with maximal ideal $\fm$ admitting a non-zero element $a\in\fm$ for which  the ideal $(0:a)$ is isomorphic to $R/aR$.
  We study minimal free resolutions of finitely generated $R$-modules $M$, with particular attention to the case when $\fm^4=0$. Let $e$ denote the minimal number of generators of $\fm$. If $R$ is Gorenstein with $\fm^4=0$ and $e\ge 3$, we show that $\Poi MRt$ is rational with denominator $\HH R{-t} =1-et+et^2-t^3$, for each finitely generated  $R$-module $M$. In particular, this conclusion applies to generic Gorenstein algebras of socle degree $3$. 
\end{abstract}

\author[I.~B.~Henriques]{In\^es~B.~Henriques}
\address{In\^es BonachoDosAnjos Henriques\\ Department of Mathematics\\
   University of Nebraska\\ Lincoln\\ NE 68588\\ U.S.A.\\ Tel.: +1-402-8057140\\
              Fax: +1-402-472-8466}
              \email{ihenriques@math.unl.edu}

\author[L.~M.~\c{S}ega]{Liana M.~\c{S}ega}
\address{Liana M.~\c{S}ega\\ Department of Mathematics and Statistics\\
   University of Missouri\\ \linebreak Kansas City\\ MO 64110\\ U.S.A.}
     \email{segal@umkc.edu}

\subjclass[2000]{Primary 13D02. Secondary 13D07}
\thanks {Research partly supported by NSF grant DMS 0803082 (I.B.H.)}
\thanks{Part of this research was done during the KRA workshop at the
University of Lincoln-Nebraska, Summer 2008, with the support of the
Nebraska MCTP grant NSF DMS-0354281 (L.M.S.)}

\maketitle

\section*{Introduction}\label{intro}
Let $R$ be a commutative noetherian local ring with maximal ideal $\fm$
and residue field $k$, and let $M$ be a finite (meaning finitely
generated) $R$-module. We study the minimal free resolution of $M$ over
$R$ by means of the {\it Poincar\'e series} of $M$, which is the formal
power series 
$$
\Poi MRt=\sum_{n=0}^{\infty}\rank_k\Tor nRMkt^n\,.
$$
The {\it Hilbert series} of $M$ is the formal power series
\[\HH M t=\sum_{n=0}^{\infty}\rank_{k}\Big(\frac{\fm ^n M}{\fm
^{n+1}M}\Big)\ t^n.\]

In this paper we  identify a class of local rings $R$ with $\fm^4=0$ with the property that  $\HH{R}{-t}\Poi M R t$ is a polynomial in $\mathbb Z[t]$ for every finite $R$-module $M$; in particular, $\Poi M R t$ is a rational function in $t$.

Research on the rationality of Poincar\'e series has uncovered the
existence of both \textit{bad} and \textit{good} behavior. Anick
\cite{Anick} constructed the first example of a local ring $R$ for which
$\Poi kRt$ is  irrational. Jacobsson \cite{Jac} produced a
local ring $R$ and finite $R$-module $M$ such that  $\Poi kRt$ is
rational, but $\Poi MRt$ is irrational. More recently, Roos \cite{Roos}
constructed a local ring $R$ and a family of finite $R$-modules
$\{M_i\}_{i\ge 0}$ such that the series $\Poi {M_i}Rt$ is  rational for
all $i\ge 0$, but do not admit a common denominator. Remarkably, in all
these examples the rings are short, in the sense that they satisfy
$\fm^3=0$. On the other hand, there is a growing gallery of classes of
rings $R$ which are {\it good}, in the sense of \cite{Roos}, that is,
there exists  $d_R(t)\in\mathbb Z[t]$ such that $d_R(t)\Poi M R t\in
\mathbb{Z}[t]$ for all finite $R$-modules $M$. Among such rings are
those with $\mu(\fm)-\depth R\le 3$ and the Gorenstein rings $R$ with
$\mu(\fm)-\depth R=4$, cf.\! \cite{AKM} and \cite{JKM}, where  $\mu(\fm)$ denotes the minimal number of generators of $\fm$.

All Gorenstein rings $(R,\fm,k)$ with $\fm^3=0$ are known to be good,
due to \!Sj\"odin \cite{SJO79}; see also Avramov, Iyengar, \c Sega
\cite{AIS}. However, Gorenstein local rings with $\fm^4=0$ can exhibit bad behavior, as  there exist examples with $\mu(\fm)\ge 12$ for which $\Poi kRt$ is  irrational (cf. \!B\o gvad \cite{Bogvad}). 

We state below our main result, in the primary case of interest.   
\begin{Theorem}
Let $(R,\fm,k)$ be Gorenstein with $\fm^4=0$ and $\mu(\fm)=e\geq3$.

If there exists a non-zero element  $a\in \fm$ such that the ideal
$(0:_{R}a)$ is principal, then the following
hold:\begin{enumerate}[\quad\rm(1)]
\item $\HH{R}{t}=1+et+et^2+t^3$;
\item $\HH{R}{-t}\Poi k R t=1$;
\item
$\HH{R}{-t}\Poi M R t$ is in $\mathbb{Z}[t]$, for every finite $R$-module $M$.
\end{enumerate}
\end{Theorem}

The hypothesis of our  theorem is verified by `generic' Gorenstein rings with $\fm^4=0$ and $\mu(\fm)\ge 3$. 
 Indeed, the set of Gorenstein standard graded $k$-algebras $R$ with $\ds{\HH R t=1+et+et^2+t^3}$ can be parametrized by a projective space by means of the inverse systems correspondence of Macaulay, which  associates to each such ring
a form of degree $3$.  With this parametrization, Conca, Rossi and Valla \cite{CRV} show that a generic Gorenstein graded $k$-algebra $R$ as above contains a non-zero element $a\in\fm$ for which the ideal $(0:_{R}a)$ is principal. They also establish $(\rm 2)$, which is equivalent to saying that $R$ is a Koszul algebra.

While the most significant statement of our results is given  for Gorenstein rings, in the body of the paper we work with larger classes of rings. We say that an element  $a\in R$ is an {\it exact zero divisor} if $R\neq(0:_{R} a)\cong R/aR\neq0$.
 Note that an element $a$ satisfying the hypothesis of the Theorem is an exact zero divisor.

The paper is organized as follows. 
In sections 1 and 2 we show that for an exact zero divisor $a$ some homological properties of $R$ and $R/aR$ are closely linked, and so are the Hilbert series of these rings.

In section 3, we establish our main result, Theorem \ref{1ee1:Koszul:Syzygy:Good}.  The rings considered are proved to be Koszul homomorphic images of a codimension two complete intersection ring under a Golod homomorphism. Results of Herzog and Iyengar \cite{HI} then yield that every finite module over $R$ has a Koszul syzygy module. (See \ref{Koszul:discussion} for the definition of a Koszul module.) In the last section we study Gorenstein rings and deduce the theorem above.

\section{Exact zero divisors and Poincar\'e series}

In this section we introduce the notion of exact zero divisor and study some 
of its properties. The main results are Theorems \ref{Poin:R:R/aR} and 
\ref{Poin:R:R/aR:Graded}.  Wide-ranging generalizations of the results of 
this section are obtained in \cite{AIS}.

We start with a remark concerning ideals with cyclic annihilators.  A version
of the equivalence of $(1)$ and $(3)$, which does not require commutativity,
is noted by Nicholson \cite{Nic04}, who calls an element $a$ satisfying these
conditions \textit{morphic}.

\begin{remark}
\label{Fcomplex}
For a commutative ring $R$ and for $a\in R$ the following are equivalent:
\begin{enumerate}
\item $(0:a)\cong R/aR$; 
\item there exists $b\in R$ such that the sequence of homomorphisms
\begin{equation*}
\ds{\F{a,b}{R}=\cdots\lra R \xra{b}R\xra{a}R \xra{b}R \xra{a}R\rightarrow 0}
\end{equation*}
is a free resolution of $R/aR$ over $R$;
\item there exists $b\in R$ such that $(0:_{R}a)=bR$ and $(0:_{R}b)=aR$
\end{enumerate}

Indeed, (1) implies $\ds{(0:_{R}a)}=bR$ for some $b\in R$, and
$bR\cong R/aR$ yields $\ds{(0:_{R}b)=aR}$, so (3) holds.
The implications (3)$\implies$(2)$\implies$(1) are evident.
\end{remark}

\begin{remark}
\label{EZD in R:principal colon ideal}
When $(R,\fm,k)$ is a local Noetherian ring, the equivalence of the 
conditions above for $a\in \fm$ is also proved by Soto.  In addition,
he proves that they are also equivalent to the freeness of the 
$R/aR$-module $(0:_{R}a)$; see \cite[Prop.~1]{Soto}.
\end{remark}

\begin{definition}\label{EPZD:rings}
Let $R$ be a commutative ring.  An element $a\in R$ is said to be an 
\textit{exact zero divisor} if it satisfies
\begin{equation}
\label{eq:epzd}
R\neq(0:_{R}a)\cong R/aR\neq0.
\end{equation}
We call any element $b$ provided by Remark \ref{Fcomplex}
a \textit{complementary divisor} for $a$,  and say that $a,b$ form an \textit{exact pair 
of zero divisors}. Note that if $R$ is local then $b$ is unique up to a unit. 
\end{definition}

\begin{bfchunk}{Exact zero divisors under flat homomorphisms.}
\label{Inflation:EPZD} 
Let $\varphi\col R\to  S$ be  a flat ring homomorphism.  It follows from the discussion above that the images under $\varphi$ of an exact pair of zero divisors in $R$ form an exact pair of zero divisors in $S$ provided they are non-zero.
\end{bfchunk}

\begin{example}
\label{exactpair:ci:example}
Let $(P,\fp,k)$ be a local ring  and $u,v,w$ be elements in $\fp$.  If $uv,w$  is a regular sequence, then the images of $u,v$ in $\ds{Q=P/(uv,w)}$ form an exact pair of zero divisors, as seen by direct computation. 
\end{example}

\begin{lemma}\label{divisorRees}
Let $R$ be a commutative  ring and $M, N$ be $R$-modules. If $a\in R$ satisfies $aM=0$ and $\ds{N/aN\cong (0:_{N}a)}$, then there exists an isomorphism of $R$-modules:
\[\Hom {R/aR} M {N/aN} \cong \Hom R M N. \]

In particular, if $R$ is local with maximal ideal $\fm$ and  $a$ is an exact zero divisor, then 
\[
\left(0:_{\ov{R}}\ov{\fm}\right)\cong(0:_{R}\fm)
\] 
where $\ov{R}=R/aR$ and $\ov{\fm}=\fm/aR.$
\end{lemma}

\begin{proof}
The following $R$-module isomorphisms are $R/aR$-linear: 
\[ N/aN\cong (0:_{N}a)\cong \Hom R{R/aR}N\,.\]
We have thus: 
\begin{align*}
\Hom{R/aR}{M}{N/aN}&\cong \Hom{R/aR}{M}{\Hom R{R/aR}{N}}\\
&\cong \Hom R {M\otimes_{R/aR}R/aR}{N}\\
&\cong \Hom R {M}{N}
\end{align*}
where the second row is given by Hom-Tensor adjointness. 

The last conclusion is obtained by taking $M=R/\fm$ and $N=R$. 
\end{proof}
\begin{remark}
\label{R/a-Gorenstein}
Lemma \ref{divisorRees} implies that if the ring $R$ is artinian and $a$ is an exact zero divisor, then $R$ is Gorenstein if and only if $\ov R=R/aR$ is Gorenstein. This statement holds more generally without the assumption that $R$ is artinian, see \cite[Thm 1]{Soto}.
\end{remark}

\begin{theorem}\label{Poin:R:R/aR}
Let $(R,\fm,k)$ be a local ring and $a,b$ be elements in $\fm\smallsetminus \fm^2$ forming an exact pair of zero divisors. For every finite $R/aR$-module $M$ one has \begin{equation}\label{large}
\Poi M R t=\Poi M{R/aR}t\cdot \frac{1}{1-t}
\end{equation}
\end{theorem}

\begin{proof}
Let $\F{a,b}{R}$ be as in Remark \ref{Fcomplex} and ${\bf G}
$ be a minimal free resolution of $k$ over $R$. The canonical map $\rho\col R/aR\to k$ extends to a morphism of complexes ${\boldsymbol{\varphi}}\col \F{a,b}{R}\to {\bf G}$, 
making the following diagram commute:
$$\xymatrix{
{\F{a,b}{R}}\ar[d]_{{\boldsymbol{\varphi}}}\ar@{->>}[r]^{\simeq}
&R/a\ar[d]^{\rho}\\
{\bf G}\ar@{->>}[r]^{\simeq} &k}$$
We show that $\varphi_{n}\otimes_{R}k$ is injective by induction on $n$ where $\varphi_{n}$ is the $n^{\text{th}}$ component of ${\boldsymbol{\varphi}}$. Indeed, if $\varphi_{2n+1}(r)\in\fm G_{2n+1}$ for some $r\in R$ then by minimality of ${\bf G}$,\ $\ds{\partial_{2n+1}(\varphi_{2n+1}(r))=a\cdot \varphi_{2n}(r)\in \fm ^2G_{2n}}$. 

\begin{equation*}\label{homomorphism}
\begin{gathered}
\xymatrix{
R \ar[d]_{\varphi_{2n+2}}\ar[r]^{b}&R \ar[d]_{\varphi_{2n+1}}\ar[r]^{a}&R\ar[d]_{\varphi_{2n}}\\
G_{2n+2}\ar[r]_{\partial_{2n+2}}&G_{2n+1}\ar[r]_{\partial_{2n+1}}&G_{2n}}
\end{gathered}
\end{equation*}
Since $a$ is part of a minimal generating set for $\fm$ we conclude that $\varphi_{2n}(r)\in\fm G_{2n}$.  By the induction hypothesis the map $\varphi_{2n}\otimes_{R}k$ is injective, hence $r$ is in $\fm$. Thus, $\varphi_{2n+1}\otimes_{R}k$ is injective. A similar argument asserts that $\varphi_{2n+2}\otimes_{R}k$ is injective, concluding the induction step. The base of induction is clear.  In this way we have that the map \[\rho_{*}\col\Tor *R{R/aR}k\lra\Tor *R k k\quad \text{ is injective.}\]
Thus, as established by Levin in \cite[1.1, (3)$\implies$(2)]{Lev80},
\[\Poi M R t=\Poi M{R/aR}t\cdot\Poi {R/aR}R t.\]  
To conclude the proof note that $\Poi {R/aR}R t=(1-t)^{-1}$\!, since $\F{a,b}{R}$ is a minimal free resolution of $R/aR$ over $R$. 
\end{proof}

We establish next a graded version of Theorem \ref{Poin:R:R/aR}. To do this we first fix some standard terminology for graded rings and modules.

\begin{bfchunk} {Graded Poincar\'e series.}
Let $k$ be a field and let ${\ds R=\oplus_{n\geq0}R_n}$ be a \textit{standard graded $k$-algebra} (in the sense of \cite[1.8]{HI}). 
 Recall that each graded module $\ds{M}$ over $\ds{R}$ has a graded free resolution and any such resolution induces a natural grading on \[\Tor i R M k=\bigoplus_{j}\Tor i R M k _{j}.\] One defines the\textit{ graded Betti numbers of $M$ over $R$} to be \[\beta_{i,j}^{R}(M)=\rank_{k}\ \Tor i R M k _{j}\]
 and the \textit{graded Poincar\'e series of $M$ over $R$} is defined to be \ \ \begin{equation*}
\Poig MRst=\sum_{i,j}\beta_{i,j}^{R}(M)\ s^i t^j\ 
\end{equation*}
 Thus $\Poig MRst$ is a formal power series in $t$ with coefficients in the Laurent polynomial ring in $s$: \[\Poig MRst\in\ds{ \mathbb{Z}\left[s^{\pm1}\right][\![t]\!]}\]

As in \cite[\textsection 2]{Sega}, we say that a finite graded module $M$ over a standard graded $k$-algebra $R$ has a \textit{linear resolution} if 
${\Tor i {R} {M} k}_j=0$ for $j\neq i$; 
equivalently,  if the graded Poincar\'e series
$\Poig M R s t$  can be written as a formal power series in the product $st$.
\end{bfchunk}
\begin{theorem}\label{Poin:R:R/aR:Graded}
Let $R$ be a standard graded $k$-algebra and $a,b$ be elements in $R_{1}$ forming an exact pair of zero divisors~on $R$. For every finite  positively graded $R/aR$-module $M$, one then has \begin{equation}
\Poig M R s t=\Poig M{R/aR} s t\cdot \frac{1}{1-st}
\end{equation}
\end{theorem}

\begin{proof}
The proof of this theorem is a graded version of the previous proof.
A book-keeping exercise shows that the proof in \cite[1.1, (3)$\Rightarrow$(2)]{Lev80} translates smoothly to the graded case. Namely, that argument uses spectral sequences that preserve grading and lead to the following conclusion:
\[\ds{\Poig{M}{R}s t=\Poig{M}{R/aR}s t\cdot\Poig{R/aR}{R}s t=\Poig{M}{R/aR}s t\cdot\frac{1}{1-st}}\ .\qedhere\]
\end{proof}

We discuss next the Koszul property for rings and modules. We refer the reader to \cite{HI} and \cite{Sega} for detailed studies of the Koszul property.

Throughout we adopt the following notation for the associated graded objects over a local ring $(R,\fm,k)$: let $\ds{\agr R}$ denote the associated graded ring and let $\ds{\agr M}$ denote the associated graded module of an $R$-module $M$; that is,
\[\ds{\agr R=\bigoplus_{n\geq 0}\ \fm ^n /\fm ^{n+1}}
\quad\text{and}\quad
\ds{\agr M=\bigoplus_{n\geq 0}\ \fm ^n M/\fm ^{n+1} M}.\]

 \begin{bfchunk}{Koszul modules.}
\label{Koszul:discussion} 
 As in \cite{HI}, we say that a finite module $M$ over a local ring $(R,\fm, k)$ is  \textit{Koszul} if for a minimal free resolution $\textbf{F}$ of $M$ the complex 
\[\operatorname{lin}^{R}(F)=\ \cdots\lra\agr{F_{n+1}}(-n-1)\lra \agr{F_{n}}(-n)\lra \cdots\lra \agr{F_{0}}\lra \agr{M}\lra 0\]
with differentials induced from $\textbf{F}$ is acyclic. A local ring $R$ is said to be \textit{Koszul} if its  residue field $k$ is Koszul as a module over $R$.

From \cite[2.3]{Sega} and \cite[1.5]{HI} the following statements are equivalent for a finite module $M$ over a local ring $R$:
\begin{enumerate}
\item
$M$ is a Koszul $R$-module
\item
$\agr M$ has a linear resolution over $\agr R$
\item
${\Tor i {\agr R} {\agr M} k}_j=0$ for $j\neq i$
\item
$\Poig {\agr M} {\agr R} s t$  can be written as a formal power series in the product $st$.
\end{enumerate}

Further, as proved in \cite[Prop. 1.8]{HI}, the Poincar\'e series of a Koszul module $M$ over a local ring $R$ satisfies
\begin{equation}\label{Poi:Koszul}
\Poi M R t=\Poi {\agr M}{\agr R}t=\HH M{-t}\HH R{-t}^{-1}.
\end{equation}

A standard graded $k$-algebra $R$ is said to be a \textit{Koszul algebra} if $k$ admits a linear free resolution as a graded module over $R$. It follows that a local ring $R$ is Koszul if and only if $\agr R$  is Koszul (as a standard graded $k$-algebra). 
 Consequently, a local ring is Koszul if and only if $\Poi {\agr k}{\agr R}{1,t}=\HH R{-t}^{-1}$\!,  as shown in \cite[Thm 1]{Fro99}.
\end{bfchunk}

\begin{bfchunk}{Graded Koszul modules.}
A graded module $M$ over a standard graded $k$-algebra $R $ is said to be Koszul if for a minimal graded free resolution $\textbf{F}$ of $M$, letting $\fm$ denote the irrelevant  maximal ideal of $R$, the complex of $R_{\fm}$-modules $\operatorname{lin}^{R}(F_{\fm})$  is acyclic. If a graded $R$-module $M$  has a linear resolution, then it is Koszul, but the converse fails in general; see \cite[1.9]{HI}. 
\end{bfchunk}

\begin{corollary}\label{1efs:Koszul:R:R/aR:Graded}
Let $R$ be a standard graded $k$-algebra and $a,b$ be elements in $R_{1}$ forming an exact pair of zero divisors on $R$. A finite positively graded $R/aR$-module $M$ has a linear resolution over $R$ if and only if it has a linear resolution over $R/aR$.
In particular, the ring $R$ is Koszul if and only if the ring $R/aR$ is Koszul.
\end{corollary}

\begin{proof} Set $\ov{R}=R/aR$.
As a consequence of Theorem \ref{Poin:R:R/aR:Graded}, $\ds{\Poig M R s t}$ is a 
power series in the product $st$ precisely when $\ds{\Poig M {\ov{R}} s t}$ is a power series in the product $st$. 
In this way, it becomes clear that $M$ has a linear resolution over $R$ if and only if it has one over $\ov{R}$. It follows immediately from the case $M=k$, that the ring $R$ is Koszul if and only if the ring $R/aR$ is Koszul.
\end{proof}

\section{Short local rings with balanced Hilbert series}

In this section we study properties of exact pairs of zero divisors over local rings $R$ satisfying $\fm^4=0$ and $H_R(-1)=0$. Theorem \ref{1efs:thm1} establishes that the initial forms $a^*,b^*$ of an exact pair of zero divisors $a,b$ in $R$ remain an exact pair of zero divisors in $\agr R$, yielding a numerical relation between the Hilbert series of $R$ and $R/aR$.

We consider both local and graded rings, hence we proceed to clarify terminology.  If $R$ is a standard graded $k$-algebra and $M=\oplus_{n\ge s}M_n$ is a finite $R$-module, then the Hilbert series of $M$ is the formal power series 
\[
\HH Mt=\sum_{n\ge s}\rank_k(M_n)t^n\,.
\]
If $R$ is artinian, then $R$ is in particular a local ring and this notion of Hilbert series coincides with the one defined for a local ring, whenever the module is generated in  degree $0$. In this setting, we will not distinguish between the local notion and the graded notion.

We say that a module $M$ has {\it balanced Hilbert series} if  $\HH M{-1}=0$.

\begin{lemma}
\label{Hilbert}
Let $R$ be a standard graded $k$-algebra and $a\in R_1$. 

If $M$ is a finite graded $R$-module such that there exists an isomorphism of graded $R$-modules 
$\ds{M/aM(-1)\cong (0:_M a)}$, then  
\begin{equation}\label{HH}
\HH Mt=(1+t)\HH {M/aM}t.
\end{equation}
In particular, if $a,b\in R_1$ form an exact pair of zero divisors then {\rm\eqref{HH}} holds with $M=R$, hence $R$ has balanced Hilbert series. 
\end{lemma}

\begin{proof}
The additivity property of the Hilbert series in the exact sequence
\begin{equation*}
0\to (0:_Ma)(-1)\to  M(-1)\xrightarrow{a}M\to M/aM\to 0
\end{equation*}
gives an equality
$$
t\HH{(0:_Ma)}t-t\HH Mt+\HH Mt-\HH {M/aM}t=0\,.
$$
From the isomorphism $\ds{M/aM(-1)\cong (0:_M a)}$ of graded $R$-modules we then have an equality
$$
(1-t)\HH Mt=(1-t^2)\HH {M/aM}t
$$
which establishes \eqref{HH}

If $a,b\in R_1$ form an exact pair of zero divisors, the complex $\F {a,b}R$ is exact and gives an isomorphism $\ds{R/aR(-1)\cong (0:_R a)}$ of graded $R$-modules yielding the desired conclusion. \qedhere
\end{proof}

Let  $(R,\fm,k)$ be a local ring. 
When $M$ is a finite  $R$-module we let $\mu(M)$ denote the minimal number of generators of $M$ over $R$; thus, one has $\mu(M)=\rank_kM/\fm M$. Also, $\lambda(M)$ denotes the length of $M$. When $\fm M=0$, one has $\mu(M)=\lambda(M)$. 

We collect below several basic facts used throughout.

\begin{remark}
\label{basic}
Let $(R,\fm,k)$ be a local ring and $a\in\fm$. 
The following then hold: 
\begin{enumerate}
\item If $a\notin\fm^2$, then $aR\cap \fm^2=a\fm$.
\item If $a\neq 0$, then $\lambda(aR)=\lambda(a\fm^2)+\mu(a\fm)+1$.
\item $\lambda(0:_Ra)=\lambda(R)-\lambda(aR)$.
\item If $\fm^4=0$, then $\lambda(R)=1+\mu(\fm)+\mu(\fm^2)+\mu(\fm^3)$.
\end{enumerate}
Indeed, to prove (1), let $x\in R$ such that $ax\in\fm^2$. Since $a$ is not contained in $\fm^2$\!, we conclude that $x$ is not an unit, hence $ax\in a\fm$. Formula (2) is given by a length count in the exact sequences: 
\begin{align*}
&0\to a\fm^2\to a\fm\to a\fm/a\fm^2\to 0\\
&0\to a\fm \to aR\to aR/a\fm\to 0\, .
\end{align*}
Formula (3) is given by a length count in the exact sequence
$$
0\rightarrow(0:_Ra)\rightarrow R\rightarrow aR \rightarrow 0\, ,
$$
and (4) by a length count on the quotients of the $\fm$-adic filtration: \[0=\fm^4\subset\fm^3\subset\fm^2\subset\fm\subset R\,.\]
\end{remark}

\begin{proposition}\label{1efs:EPZD}
Let $(R,\fm,k)$ be a local ring with $\fm^4=0$ and  balanced Hilbert series.
 For $a,b$ in $\fm\smallsetminus \fm^2$, the following statements are equivalent:
\begin{enumerate}[\quad\rm(1)]
\item $a,b$ is an exact pair of zero-divisors
\item $(0\col a)=bR$
\item $ab=0$, $a\fm^2=\fm^3=b\fm^2$ and $\mu(a\fm)=\mu(\fm)-1=\mu(b\fm)$.
\end{enumerate}
\end{proposition}

\begin{proof}
It is clear that (1) implies (2). For the rest it suffices to show that (2) and (3) are equivalent; indeed, if this is known, we see that (2) implies $(0\col b)=aR$, since (3) is symmetric with respect to $a$ and $b$.

We may assume $ab=0$ as this holds under any of the conditions $(\rm 1)$ to $(\rm 3)$.

Set $e=\mu(\fm)$ and $s=\mu(\fm^3)$. Since  $\fm^4=0$, one gets $\mu(a\fm^2)=\lambda(a\fm^2)$ and $\mu(b\fm^2)=\lambda(b\fm^2)$. Remark \ref{basic}(2) gives
\begin{equation}
\label{eq1length}
\lambda(aR)=\mu(a\fm^2)+\mu(a\fm)+1\,.
\end{equation}
Since $b\not\in\fm^2$, we may choose $b, x_2, x_3,\dots, x_e$ to be a minimal system of generators for $\fm$. As $ab=0$, we see that $a\fm$ is generated by $ax_2, \dots, ax_e$, hence 
\begin{equation}
\label{ineq1}
\mu(a\fm)\le e-1\,.
\end{equation}
From the inclusion $a\fm^2\subseteq \fm^3$ we obtain the obvious inequality
\begin{equation}
\label{ineq2}
\mu(a\fm^2)\le s
\end{equation}
with equality if and only if $a\fm^2=\fm^3$. 
Using \eqref{ineq1} and \eqref{ineq2}  in \eqref{eq1length} we thus have
\begin{equation}
\label{ineq3}
\lambda(aR)\le e+s
\end{equation}
with equality if and only if $\mu(a\fm)= e-1$ and $a\fm^2=\fm^3$. Similarly, we obtain
\begin{equation}
\label{ineq4}
\lambda(bR)\le e+s
\end{equation}
with equality  if and only if $\mu(b\fm)= e-1$ and $b\fm^2=\fm^3$.
The hypothesis $\ds{\HH R{-1}=0}$ implies $\lambda(R)=2(e+s)$ and Remark \ref{basic}(3) gives 
 \[\lambda(0:a)=2(e+s)-\lambda(aR)\,.\] 
Thus  the inequality \eqref{ineq3} yields
\begin{equation}
\label{ineq5}
\lambda(0:a)\ge e+s
\end{equation} 
with equality if and only if $\mu(a\fm)= e-1$ and $a\fm^2=\fm^3$.

Since $ab=0$, we have $bR\subseteq (0:a)$ and the inequalities above yield
$$
\lambda(bR)\le e+s\le \lambda(0:a).
$$
We have thus $bR=(0:a)$  if and only if equalities hold in \eqref{ineq4} and \eqref{ineq5}, that is, if and only if condition (3) is satisfied.  
\end{proof}

\begin{lemma}
\label{EZD:CZD:m2}
Let $(R,\fm,k)$ be a local ring satisfying $\fm^4=0$ and $\mu(\fm)=e\ge 3$. 
\begin{enumerate}[\quad\rm(1)]
	\item If $R$ has balanced Hilbert series and an exact zero divisor lies in $\fm\smallsetminus\fm^2$ then so do its complementary divisors.
	\item If $\mu(\fm^3)+2\le e$, then a non-zero exact zero divisor is not contained in $\fm^2$.
\end{enumerate}
\end{lemma}

\begin{proof} Let $f=\mu(\fm^2)$ and $g=\mu(\fm^3)$.

(1) As $\fm^4=0$, we have $\HH R{t}=1+et+ft^2+gt^3$. If $\HH R{-1}=0$ then  $g=f-e+1$. Now let $a\in\fm\smallsetminus\fm^2$ be an exact zero divisor and $b$ its complementary divisor. 

Assume that $b$ is in $\fm^2$. One then has 
 \[aR=(0:_{R}b)\supseteq(0:_{R}\fm^2)\supseteq\fm^2.\]
 Hence $\fm^2=a\fm$ by Remark \ref{basic} (1) and 
\[
f=\mu(\fm^2)=\mu(a\fm)\le \mu(\fm)=e.
\]
We have thus    $g=f-e+1\le 1$. Since $e\ge 3$, a contradiction is derived from (2).

(2)  Assume $b$ is a non-zero exact zero divisor in $\fm^2$ and let $a$ be a complementary zero divisor. Then $b\fm\subseteq\fm^3$ hence 
\[
\lambda(bR)=\lambda(bR/b\fm)+\lambda(b\fm)\leq \mu(bR)+\lambda(\fm^3)=1+g.
\]
On the other hand, $a\fm\subseteq\fm^2$ as $a\in\fm$, so 
\[
\lambda(aR)=\lambda(aR/a\fm)+\lambda(a\fm)\leq \mu(aR)+\lambda(\fm^2)= 1+\mu(\fm^2)+\mu(\fm^3)=1+f+g.
\]

From Remark \ref{basic} (4) and (3), we have
\[1+e+f+g=\lambda(R)=\lambda(aR)+\lambda(bR)\leq 1+f+g+1+g
\]
implying $e\leq g+1$, a contradiction. 
\end{proof}

We use the notation introduced in Section 1: $\agr R$ and $\agr M$ denote the associated graded ring and module respectively.

  Given an ideal $I$ of $R,$ let $I^*$ denote the ideal of $\agr R$ generated by the initial forms of elements in $I$. Also, let $r^*$ denote the initial form of an element $r$ in $R$.

\begin{theorem}\label{1efs:thm1} 
Let $(R,\fm,k)$ be a local ring with $\fm^4=0$ and balanced Hilbert series. 
If $a, b\in\fm\smallsetminus \fm^2$ form an exact pair of zero divisors then the following hold:
\begin{enumerate}[\quad\rm(1)]
\item their initial forms $\ds{a^*,b^*}$ in $\ds{\agr R}$ form an exact pair of zero divisors;
\item the graded rings $\ds{\agr {(R/aR)}}$ and $\ds{\agr R/a^*\agr R}$ are naturally isomorphic;
\item
$\ds{\HH{\agr R /a^*\agr R}t
=\HH{R/aR}t=\HH{R}{t}(1+t)^{-1}}$.
\end{enumerate}
\end{theorem}

\begin{proof} (1) Note that $(\agr R,\agr \fm)$ is a local ring, where  $\agr\fm$ denotes  the irrelevant maximal ideal $\agr[\geqslant 1]R$ of  $\agr R$. One has  
\[\mu(a\fm)=\lambda(a\fm /a\fm ^2)=\mu(a^*\agr{\fm})\quad\text{and}\quad \mu(\fm)=\lambda(\fm /\fm ^2)=\mu(\agr{\fm})\,.\] 
Further, the property $a\fm^2=\fm ^3=b\fm^2$ given by Proposition \ref{1efs:EPZD} implies
 \[ a^*\:{(\agr{\fm})}^2=({\agr{\fm}})^{3}=b^*\:{(\agr{\fm})}^2.\]
Proposition \ref{1efs:EPZD} shows that $a^*,b^*$ form an exact pair of zero divisors in $\agr R$.

(2) The canonical homomorphism $R\to R/aR$ induces a surjective homomorphism of graded $k$-algebras
$\varphi\col \agr R \to \agr{\left(R/aR\right)}$.  From \cite[(1.1)]{VaVa},  $\ds{\ker \varphi=a^*\agr R}$ if and only if \begin{equation}\label{vv}
\fm ^{i}\cap aR=a\fm ^{i-1}\quad\text{for all}\ \ i\geq 1\,.
\end{equation}

The cases $\ds{i=1}$ and $\ds{i\geq4}$ are obvious and case $i=2$ is given by Remark \ref{basic}(1). Finally, $i=3$ follows from Proposition \ref{1efs:EPZD}(3) and we conclude that \eqref{vv} holds.

(3) Note that $\HH {\agr R}t=\HH Rt$ and $\HH {{\agr R}/a^*{\agr R}}t=\HH {\agr{(R/aR)}}t=\HH {R/aR}t$. The statement follows then from (1) and Lemma \ref{Hilbert}.  
\end{proof}

\section{Modules over short local rings}
\setcounter{equation}{0}
In this section we study the Koszul property of finite modules over local artinian rings $(R,\fm,k)$ with $\fm^4=0$ and balanced Hilbert series.  

We establish our main result, Theorem \ref{1ee1:Koszul:Syzygy:Good}, whose primary case of interest (over Gorenstein rings) is discussed in section 4.
 An important structural result is Proposition \ref{General:Golodmap} yielding that, in the presence of special elements $a,b,c\in R$, the rings considered are homomorphic
images of a codimension $2$ complete intersection by means of a Golod homomorphism. This allows to apply results of Herzog and Iyengar \cite[\S 5]{HI} on the Koszul property. 

\begin{proposition}\label{1efs:Koszul:R:R/aR}
Let $(R,\fm,k)$ be a local ring with $\fm^4=0$ and balanced Hilbert series.  If $a,b$ in $\fm\smallsetminus \fm^2$ form an exact pair of zero divisors and $M$ is an $R/aR$-module, then
$M$ is Koszul over $R$ if and only if it is Koszul over $R/aR$.

In particular, the ring $R$ is Koszul if and only if the ring $R/aR$ has this property.
\end{proposition}

\begin{proof}
Recall from \ref{Koszul:discussion} that a module is Koszul over a local ring if and only if 
its associated graded module has a linear resolution over the associated graded ring.\ 
From Theorem \ref{1efs:thm1}(1), the initial forms  $a^*,b^*$ form an exact pair of zero divisors in $\agr R$.  As $\agr M$ is a positively graded $\agr{R}/a^*\agr R$-module, Corollary \ref{1efs:Koszul:R:R/aR:Graded} shows that $\agr M$ has a linear resolution over $\agr R$ precisely when $\agr M$ has a linear resolution over $\agr R/a^*\agr R$. This is the desired conclusion because the graded rings $\agr R/a^*\agr R$ and $\agr{\left(R/aR\right)}$ are naturally isomorphic by Theorem \ref{1efs:thm1}(2).
\end{proof}

Recall that the $n^{\textit{th}}$ \textit{syzygy module} of a finite $R$-module $M$ is the cokernel of the $n^{\text{th}}$ differential in a minimal free resolution of $M$ over $R$. It is defined uniquely up to isomorphism and we denote it $\syz nM$. If an $R$-module $N$ is isomorphic to $\syz nM$ for some $n\geq 0$, then we say that $N$ is a \textit{syzygy module of $M$}.

\begin{bfchunk}{Koszul syzygy modules.}\label{koszulsyzygymodule}  If $(R,\fm,k)$ is a zero-dimensional local ring and $M$ a finite $R$-module admitting a Koszul syzygy module  then $$\HH R {-t }\Poi M R t\in \mathbb{Z}[t].
$$ This follows immediately from \cite[Prop. 1.8]{HI}.

In this way, if $R$ is a zero-dimensional local ring over which every finite $R$-module  $M$ admits a Koszul syzygy module, then for each finite $R$-module $M$ there exits a polynomial $p_M (t)\in \mathbb{Z}[t]$ satisfying 
\[
\Poi M R t
=p_M (t)/{\HH{R}{-t}}.
\]

In general an uniform bound for the degree of $p_M (t)$ can not be derived. Indeed, if $R$ is a self-injective non-regular Koszul ring, then for every $n\ge 0$ there exists a non-zero finite $R$-module $M$ satisfying $k\cong \Syz R n M$. Since $\Poi k R t=\HH{R}{-t}^{-1}$ one then has
\[
\Poi M R t
=t^n\Poi k R t+q_M (t)
=\frac{t^n}{\HH{R}{-t}}+q_M(t)
=\frac{p_M (t)}{\HH{R}{-t}}
\] 
for some (uniquely determined)  polynomials $q_M (t)$ and $p_M (t)$ of  degrees $n-1$ and $\ds{n-1+\deg \HH R t}$ respectively.
\end{bfchunk}

Our main result is as follows.

\begin{theorem}\label{1ee1:Koszul:Syzygy:Good}
Let $(R,\fm,k)$ be a local ring with $\HH{R}{t}=1+et+et^2+t^3$ and set $s=\rank_k(0:_R\fm)$. Assume there exists a non-zero exact zero divisor $a\in\fm$. The following statements are then equivalent:
\begin{enumerate}[\quad\rm(1)]
\item $R$ is Koszul.
\item $e\geq s+2$. 
\end{enumerate}
When they hold every finite $R$-module $M$ has a syzygy module that is Koszul.
\end{theorem}

The proof of Theorem \ref{1ee1:Koszul:Syzygy:Good} is deduced from Theorem \ref{Goodness}.

\begin{corollary}
Let $(R,\fm,k)$ be a local ring with  $\HH Rt=1+et+et^2+t^3$ and $e\ge 3$.  Assume there exists a non-zero exact zero divisor $a\in\fm$ and $\agr R$ is quadratic. 
 The ring  $R$ then is Koszul and every finite $R$-module $M$ has a syzygy module that is Koszul. 
\end{corollary}

\begin{proof}
Let $b$ be a complementary divisor for $a$, then $a,b$ form an exact pair of zero divisors in $\fm\smallsetminus\fm^2$ by Lemma \ref{EZD:CZD:m2}. Theorem \ref{1efs:thm1}(1) then establishes that the initial forms $a^*,b^*$ of $a,b$ form an exact pair of zero divisors in $\agr R$. Further, $\agr R/a^*\agr R$ is quadratic and $\HH {\agr R/a^*\agr R}t=1+(e-1)t+t^2$ by Theorem \ref{1efs:thm1}(3). By \cite[2.12]{CRV}, it follows that
$$
\rank_k(0:_{\agr R}\agr \fm)\le e-2
$$ 
Applying Theorem \ref{1ee1:Koszul:Syzygy:Good} to the local ring $\agr R$, we conclude that $\agr R$ is Koszul, hence $R$ is Koszul as well. Thus $R$ satisfies conditions (1)-(2) of the same theorem, concluding the proof.
\end{proof}

\begin{remark}
\label{connection:CRV}
Let $R$ be an artinian graded $k$-algebra with irrelevant maximal ideal $\fm$. In particular, $(R,\fm,k)$ is a local ring. In \cite{CRV}, Conca, Rossi and Valla consider standard graded $k$-algebras $R$ with  $\HH R{t}=1+et+et^2+t^3$, where $e\geq 3$, that contain elements $l,m\in R_1$ with  
\begin{equation}
\label{lm}
lm=0\quad\text{ and}\quad  \rank_k(lR_1)=\rank_k(mR_1)=e-1\,.
\end{equation}
When $R$ is quadratic or Gorenstein, Proposition \ref{1efs:EPZD} and the proof of \cite[2.13]{CRV} show  that a pair $l,m\in R_1$ satisfies \eqref{lm} if and only if it is an  exact pair of zero divisors. The Gorenstein case also follows from  Proposition \ref{Gorenstein:EZD} in Section 4. 
\end{remark}

 We introduce terminology to handle the technical ingredients of the proof.

\begin{bfchunk}{Conca generators.}\label{Conca}
Following \cite{AIS}, we say that an element $c\in \fm$ is a \textit{Conca generator modulo $aR$} if it satisfies
\[\fm ^2 +aR = c\fm+aR,\quad c\notin aR\quad\text{and}\quad c^2\in aR.\]
Let $\ov{c}$ be the image of $c$ in $\ov{R}$. These conditions imply ${\ov{\fm}}^3=0$ and $\ov{c}\notin {\ov{\fm}}^2=\ov{c}\ov{\fm}$ by Nakayama's Lemma. Further, they are equivalent to
\[\fm ^2 \subseteq c\fm+a\fm,\quad c\notin aR\quad\text{and}\quad c^2\in a\fm\] whenever $a$ is part of an exact pair of zero divisors,
as given by  \eqref{vv}.

In the terminology of \cite{AIS}, the image in $R/aR$ of a Conca generator modulo $aR$ is exactly a Conca generator of the maximal ideal $\fm/aR$.
\end{bfchunk}

\begin{bfchunk} {\it Golod homomorphisms.} The notion of Golod homomorphism was  introduced by Levin \cite{Lev75}.  An epimorphism $\psi\col(Q,\fq,k)\thra (R,\fm,k)$ is said to be {\it Golod} if the following equality holds 
\begin{equation}
\label{Golod-formula}
\ds{\Poi{k}{R}t=P_{k}^{Q}(t)(1-t(P_{R}^{Q}(t)-1))}^{-1}.
\end{equation}
\end{bfchunk}

In the proof of the next proposition we use a result of Valabrega and Valla:
\begin{chunk}{\it Initial forms and regularity.}\label{VV:cor2.7}
To establish the regularity of a sequence in a local ring $Q$ it suffices to establish regularity of the corresponding sequence of initial forms in $\agr Q$; see \cite[2.7]{VaVa}. 
\end{chunk}

As $R$ is artinian, \text{Cohen's Structure Theorem} yields a regular local ring $(P,\fp,k)$ and a surjective local homomorphism 
\[
\pi\col\ (P,\fp,k)\thra (R,\fm,k)
\]
such that $\ker \pi\subseteq \fp^2$. 

\begin{proposition}\label{General:Golodmap}
Let $(R,\fm,k)$ be a local ring 
with $\fm^4=0$ and balanced Hilbert series.  Let $a,b,c$ be elements in $\fm\smallsetminus\fm^2$ such that $a,b$ form an exact pair of zero divisors and $c$ is a Conca generator modulo $aR$.

	The ring $R$ then is Koszul, and there exist  $u,v\in \fp\smallsetminus\fp ^2$ and $w\in\fp^2$ such that  $uv,w$ is a $P$-regular sequence contained in $\ker \pi$ and such that the homomorphism \[Q=P/(uv,w)\thra R\] induced by $\pi$ is Golod.
\end{proposition}

\begin{proof}
Let $(\ov{R},\ov{\fm})$ denotes the local ring $(R/aR,\fm/aR)$.  For each $x\in R$ we let $\ov{x}$ denote the image of $x$ in $\ov{R}$. The ring  $\ov{R}$ is Koszul by \cite[1.1]{AIS}, so Proposition \ref{1efs:Koszul:R:R/aR} shows that $R$ is Koszul, and from \eqref{Poi:Koszul} we obtain  
\begin{equation}
\ds{\Poi k R t =\HH{R}{-t}^{-1}.}
\end{equation}

The rest of the proof is conducted in three steps.
The first one is to show that $\pi$ factors through a map $\psi\col Q\thra R$ where $Q$ is a complete intersection of the form $P/(uv,w)$ for some regular sequence $uv,w$ in $P$ such that $\pi(u)=a$. 

We let $(\ov P,\ov{\fp})$ denote the regular local ring $(P/uP, \fp/uP)$ and for every $x\in P$ we denote $\ov x$ the image of $x$ in $\ov P$. 

In a second step we remark that the induced map  $\ov{\pi}\col \ov P\thra \ov{R}$ is a Cohen presentation for $\ov{R}$ with $\ker \ov\pi\subseteq \ov \fp^2$ and maps a minimal generator for $\ov\fp$ to a Conca generator for $\ov\fm$. We then use \cite[1.4]{AIS} to conclude that the map $\ov{\psi}\col\ov{Q}\thra \ov{R}$ induced by $\psi$ on $\ov{Q}=Q/\tilde{u}Q$ is Golod, where $\tilde{u}$ denotes the image of $u$ in $Q$. 

The final step consists of establishing the equality \eqref{Golod-formula}
and therefore concluding that $\psi$ itself is a Golod homomorphism.

\begin{Step1}
We show there exist elements $u,v$ in $\ds{\fp\smallsetminus\fp^2}$ and $w$ in $\fp^2$ such that $\pi(u)=a,\:\pi(v)=b$ and $uv,w$\ forms a $P$-regular sequence.
The map $\pi$ then factors through a homomorphism \[\psi\col (Q,\fq,k)\thra (R,\fm,k)\ \text{where}\ \ds{Q=P/(uv,w)}.\]
\end{Step1}

Indeed, as $c$ is a Conca generator modulo $aR$, there exists $d\in \fm$ such that 
\begin{equation}\label{c2ad}
	c^2=ad\ \ \text{ and  }\ \ \fm ^2+(a)=c\fm +(a).
\end{equation}
One then has $\ov{c}\notin \ov{\fm}^2$ so $a,c$ can be extended to a minimal generating set for $\fm$. There exists thus  a minimal system of generators $x_1,\dots,x_e$ for $\fp$ satisfying $\ds{\pi(x_1)=a}$ and $\ds{\pi(x_2)=c}$. 
Set $X_i={x_i}^*$ in $\agr P$ and note that $\agr P$ is the polynomial ring over $k$ on $X_1,\dots,X_e$, since $P$ is regular.

Pick $y \in \fp$ such that $\pi(y)=d$ and set 
\[
u=x_1\quad\text{and}\quad w=x_2^2-y x_1.
\]

To carry out the first step we consider three cases. 
\begin{Case1}
$\rank_{k}\{a^{*},b^{*}\}=1$.
 
We have thus $a^{*}\agr R=b^{*}\agr R$. We may then assume $b=a+f$ for some $f\in \fm^2$, since exact pairs of zero divisors are preserved up to multiplication by a unit in $R$. Choose $z\in \fp^2$ so that $\pi(z)=f$; the map $\pi$ then factors through the map \begin{equation}
\label{Q1}
Q=\frac{P}{(x_1(x_1+z),{x_2}^2-y  x_1)}\ \xra{ \ \psi\ } R.
\end{equation}
Set $\ds{v=x_1+z}$. We claim that the sequence\ $\ds{uv,\ w}$\ is $P$-regular.
In $\ds{\agr P}$ one has \[(uv)^*={X_1}^2\quad \text{and} \quad\ds{w^*=\left\{\begin{array}{ll}
	{X_2}^2 -y^* X_1, &\text{ if }y  \in \fp\smallsetminus \fp^2\\
	{X_2}^2  ,& \text{ if }y \in \fp^2	
	\end{array}
	\right.
	}.\] 
Since both sequences \[\ds{{X_1}^2 ,{X_2}^2}\quad \text{and} \quad\ds{{X_1}^2,\ {X_2}^2-y^* X_1}\] are $k[X_1,\dots,X_e]$-regular, we conclude from \ref{VV:cor2.7} that $\ds{uv,w}$ is $P$-regular.
\end{Case1}
 \begin{Case2}$\rank_{k}\{a^{*},b^{*},c^{*}\}=3$.
 
 As $a,b,c$ can be extended to a minimal generating set for $\fm$,  then we may assume $\ds{\pi(x_3)=b}$. The map $\pi$ then factors through the homomorphism \begin{equation} 
\label{Q2}
Q=\frac{P}{(x_1 x_3,{x_2}^2-y  x_1)}\xra{\ \ \psi\ } R.
\end{equation}
Set $\ds{v=x_3}$. We claim that the sequence\ $\ds{uv,\ w}$\ is $P$-regular. In $\ds{\agr P}$ one has \[(uv)^*=X_1 X_3\quad \text{and} \quad\ds{w^*=\left\{\begin{array}{ll}
	{X_2}^2 -y^* X_1, &\text{ if }y  \in \fp\smallsetminus \fp^2\\
	{X_2}^2  ,& \text{ if }y \in \fp^2	
	\end{array}
	\right.
	}.\] 
Since both sequences \[\ds{X_1 X_3 ,{X_2}^2}\quad \text{and} \quad\ds{X_1 X_3,\ {X_2}^2-y^* X_1}\] are $k[X_1,\dots,X_e]$-regular, we conclude from \ref{VV:cor2.7} that $\ds{uv,w}$ is $P$-regular.
\end{Case2}
\begin{Case3} 
$\rank_{k}\{a^{*},b^{*},c^{*}\}=2=\rank_{k}\{a^{*},b^{*}\}$.

In other terms, $\ov{c}=\ov{b}+\ov{g}$ \ in $R/aR$ \ for some \ $g\in \fm^2$. 
Since $\ov{\fm}^3=\ov{a\fm^2}=0$, we have
\[\ov{b}^2=\ov{c}^2=0 \quad \text{and} \quad  \ov{\fm}^2=\ov{\fm}\ \ov{c}=\ov{\fm}\ \ov{b}.\]
Thus we may replace $c$ with $b$ and then \eqref{c2ad} yields $b^2=ad$.

If $d\in bR$, then $b^2=ad=0$\ which implies\ $b\in (0:_{R}b)=aR$. As $b\notin \fm^2$, it would then fall into \textit{Case 1}. It remains to treat the case when $d\notin bR$.

If $d\notin bR$, the map $\pi$ then factors through the homomorphism \begin{equation}
\label{Q3}
Q=\frac{P}{(x_1x_2,{x_2}^2-y  x_1)}\xra{\ \ \psi\ } R.
\end{equation}
Set $\ds{v=x_2}$. The sequence\ $\ds{uv,\ w}$\ is $P$-regular: because $P$ is a UFD, $x_1,x_2$ are irreducible elements and $\ds{x_2\nmid y}$.
\end{Case3}

\begin{Step2}
We observe that the map $\ov{\psi}\col \ov{Q}\thra \ov{R}$ induced by $\psi\col  Q\thra R$ \ on $\ov{Q}=Q/\tilde{u}Q$ \ is Golod.
\end{Step2}

By choice, $u=x_1$ is not in $\fp^2$, so the ring $\ov{P}=P/uP$ is regular. The induced map $\ov{\pi}\col\ov{P}\thra\ov{R}$ is a Cohen presentation with $\ker \ov\pi\subseteq \ov\fp^2$  such that the element $\ov{x}_2\in \ov{\fp}\smallsetminus \ov{\fp}^2$ maps to the Conca generator $\ov{c}$ of $\ov{\fm}$. Note that there are canonical isomorphisms 
\begin{equation}
\label{isos}
\ov{Q}=\frac{Q}{\tilde{u}Q}\cong \frac{P}{(u,w)}
\cong \frac{\ov{P}}{\ov{x}_2 ^2\ \ov{P}}\,.
\end{equation}

Under this identification the map $\ov{\pi}\col\ov{P}\thra\ov{R}$ factors through $\ov{\psi}\col \ov{Q}\thra \ov{R}$, and it follows from \cite[1.4]{AIS} that $\ov{\psi}$ is a Golod homomorphism.

\begin{Step3}
We now prove that the homomorphism $\psi\col Q\thra R$ is Golod.
\end{Step3}

We first note that, from Remark \ref{exactpair:ci:example}, the images $\ds{\tilde{u},\tilde{v}}$ of $\ds{u,v}$ in $Q$ form an exact pair of zero divisors. Thus, as in Remark \ref{Fcomplex}, 
$\ds{\F{\tilde{u},\tilde{v}}{Q}}$\ is a minimal free resolution of\ $\ov{Q}$ over $Q$. Since $\psi(\tilde{u})=a$ and $\psi(\tilde{v})=b$, one has an isomorphism of complexes 
\begin{equation*}\label{resolution2}
\F{\tilde{u},\tilde{v}}{Q} \otimes_{Q}R
\cong \F{a,b}{R},
\end{equation*}
which induces isomorphisms
\[\Tor{n} Q {\ov{Q}} R =\HH{n}{\F{\tilde{u},\tilde{v}}{Q} \otimes_{Q}R}=\HH{n}{\F {a,b} R}=0 , \ \text{ for all} \ n>0.\]
Therefore,\begin{equation}\label{Poin:Q:Q/uQ}
P_{\ov{R}}^{\ov{Q}}(t)=P_{R}^{Q}(t).
\end{equation}
Since $\ov{P}$ is a regular local ring, \eqref{isos} and \cite[Thm 6]{Tate} 
 yield equalities
\begin{align}\label{PkQ}
P_{k}^{Q}(t)&=(1+t)^{e}(1-t^2)^{-2},\\
\label{PkQbar}
P_{k}^{\ov{Q}}(t)&=(1+t)^{e-1}(1-t^2)^{-1}\,,
\end{align}
where $e=\mu(\fm)=\mu(\ov\fm)+1$. 

Now we can write the following sequence of equalities
\begin{align*}
\Poi{k}{R}t&\ds{=\frac{\Poi{k}{\ov{Q}}{t}}{(1-t(P_{\ov{R}}^{\ov{Q}}(t)-1))(1-t)}}=\frac{(1+t)^{e-1}}{(1-t(P_{\ov{R}}^{\ov{Q}}(t)-1))(1-t^2)(1-t)}\\
&\ds{=\frac{(1+t)^{e} }{(1-t(P_{R}^{Q}(t)-1))(1-t^2)^2}=\frac{P_{k}^{Q}(t)}{1-t(P_{R}^{Q}(t)-1)}}
\end{align*}
The first equality above follows from Step 2 and Theorem \ref{Poin:R:R/aR}, the rest from \eqref{PkQbar}, \eqref{Poin:Q:Q/uQ} and \eqref{PkQ} respectively.
The resulting equality yields that the map $\ds{\psi\col  Q \thra R}$ is Golod.
\end{proof}

\begin{remark}\label{inflation:Koszul}
Let $\ds{\varphi\col (R,\fm,k)\hookrightarrow(R',\fm',k')}$ be an \textit{inflation}; that is, a flat homomorphism of local rings such that $\varphi(\fm)=\fm'$. One then has $\HH{R}{t}=\HH{R'}{t}$.

For every local ring $(R,\fm,k)$ there exits a local ring $(R',\fm',k')$ and an inflation $R\hookrightarrow R'$ where $k'$ is an algebraically closed field $k$, see \cite[AC IX.41]{BourAC9}.  

An inflation $\ds{\ov{\varphi}\col  (\ov{R},\ov{\fm},k)\hookrightarrow(\ov{R}',\ov{\fm}',k')}$, where $\ov{R}=R/aR$, $\ov\fm=\fm/aR$ and $(\ov{R}',\ov{\fm}',k')$ is a local ring, lifts to an inflation $\ds{\varphi\col (R,\fm,k)\hookrightarrow(R',\fm',k')}$ satisfying  $R'/\varphi(a)R'\cong\ov{R}'$, see \cite[AC IX.40]{BourAC9}.

Given a minimal free resolution $\textbf{F}$ of an $R$-module $M$, the complex  $\ds{\textbf{F}\otimes_{R}R'}$ is a minimal free resolution of the $R'$-module $\ds{M'=M\otimes_{R}R'}$. Thus one gets 
\[\Syz{R'}{n}{M'}\cong\syz nM\otimes_{R}R'\quad\text{and}
\quad \Poi {M'}{R'}{t}=\Poi {M}{R}{t}.\] 

Further, as noted in \cite[1.8]{AIS}, the extended module $\ds{M'}$ is Koszul over $R'$ precisely when $M$ is Koszul over $R$.  
\end{remark}

\begin{theorem}\label{Goodness}
Let $(R,\fm,k)$ be a local ring with ${\fm}^4=0$ and $\HH R{-1}=0$.\\
Assume that there exist an inflation $\ds{(R,\fm,k)\to (R',\fm',k')}$ and elements $a,b,c$ in $\ds{\fm'\smallsetminus{\fm'}^2}$ such that $a,b$ form an exact pair of zero divisors in $R'$ and $c$ is a Conca generator modulo $a R'$.

 The ring $R$ then is Koszul and every finite $R$-module $M$ has a syzygy module that is Koszul.  In particular, $\Poi M Rt\HH R {-t }^{-1}\in\mathbb{Z}[t]$.
\end{theorem}

\begin{proof}
In view of Remark \ref{inflation:Koszul} and \ref{Inflation:EPZD}, we can  replace $\ds{(R',\fm',k')}$ with $\ds{(R,\fm,k)}$ and $\ds{M\otimes_{R}R'}$ with $M$.
 The ring $R$ then satisfies the hypothesis of Proposition \ref{General:Golodmap}, and hence is Koszul. 

Proposition \ref{General:Golodmap} and \cite[5.9]{HI} show that every finite \text{$R$-module} $M$ has a syzygy module that is Koszul over $R$; hence $\Poi M Rt\HH R {-t }^{-1}\in\mathbb{Z}[t]$ from \ref{koszulsyzygymodule}.
\end{proof}

\begin{proof}[Proof of Theorem {\rm \ref{1ee1:Koszul:Syzygy:Good}}]
Let $b$ be a complementary divisor for $a$. From Lemma \ref{EZD:CZD:m2}(2) $a,b$ form an exact pair of zero divisors in $\fm\smallsetminus\fm^2$, since $\mu(\fm^3)=1= 3-2\le e-2$.

Recall that $\ov{R}=R/aR$ is a local ring with maximal ideal $\ds{\ov{\fm}=\fm /aR}$ and $\HH{\ov{R}}{t}=1+(e-1)t+t^2$ as given by Theorem \ref{1efs:thm1}(3), then 
 \[\rank_{k}(\ov{\fm}/\ov{\fm}^2)=e-1\ \text{  and }\ \rank_{k}(\ov{\fm}^2)=1.\]
 Lemma \ref{divisorRees}  yields 
\[\rank_{k}(0:_{\ov{R}}\ov{\fm})=s
\] 
and \cite[4.1]{AIS} shows  that $\ov R$ is Koszul if and only if 
$$
\rank_{k}(0:_{\ov{R}}\ov{\fm})\le \rank_{k}(\ov{\fm}/\ov{\fm}^2)-1
$$
that is,  $s\le e-2$. Proposition \ref{1efs:Koszul:R:R/aR} gives that $R$ is Koszul if and only if $\ov R$ is Koszul, establishing thus the equivalence of (1) and (2). 

Assuming that (1) and (2) hold, then  \cite[Thm 4.1,(v)$\implies$(i)]{AIS} and Remark \ref{inflation:Koszul} establish the existence of an inflation $\ds{\varphi\col  (R,\fm,k)\hookrightarrow(R',\fm',k')}$ such that $\fm'\smallsetminus {\fm'}^2$ contains an element $c$ that is a Conca generator modulo $\varphi(a)R'$. Finally, from \ref{Inflation:EPZD} we can apply Theorem \ref{Goodness}
 concluding the proof.
\end{proof}

\section{Short local Gorenstein rings}

In this section we study exact pairs of zero divisors in local Gorenstein rings $(R,\fm,k)$ with $\fm^4=0$. We prove the Gorenstein version of our main result as presented in the introduction.

Recall that $R$ is Gorenstein if and only if $\ov R=R/aR$ is Gorenstein, as discussed in Remark \ref{R/a-Gorenstein}.

\begin{proposition}
\label{Gorenstein:EZD}
Let $(R,\fm,k)$ be a Gorenstein local ring satisfying $\fm^4=0$ and $\mu(\fm)=e\ge3$. Given $0\ne a\in\fm$ the following statements are equivalent:
\begin{enumerate}[\quad\rm(1)]
\item $(0:_{R}a)$ is a principal ideal;
\item $a$ is an exact zero divisor;
\item $a$ is in $\fm\smallsetminus\fm^2$ and there exists $\ds{b\in\fm\smallsetminus \fm^2}$ such that the following equivalent conditions hold:
\begin{enumerate}[\quad\rm(3.1)]
\item\ $a,b$ form an exact pair of zero divisors;
\item\ $ab=0$ and $\mu(a\fm)=e-1=\mu(b\fm)$\,.
\end{enumerate}
\end{enumerate}
If any of the conditions {\rm (1)-(3)} holds, then $\HH Rt=1+et+et^2+t^3$.
\end{proposition}´

\begin{proof}
We first show that given $a,b\in\fm\smallsetminus\fm^2,$  {\rm (3.1)} implies $\HH Rt=1+et+et^2+t^3$. Set $(\ov{R},\ov{\fm})=(R/aR, \fm/aR)$. The hypothesis that $a\in\fm\smallsetminus\fm^2$ implies
$\mu(\ov\fm)=e-1$. 
Since $R$ is Gorenstein,  Remark \ref{R/a-Gorenstein} shows  that $\ov
R$ is Gorenstein.  Note that 
\begin{equation}
\label{am}
\ov\fm^2=(\fm^2,aR)/aR\cong \fm^2/aR\cap\fm^2=\fm^2/a\fm
\end{equation}
where the equality $aR\cap\fm^2=a\fm$ is given by Remark \ref{basic}(1).

Assume $\ov\fm^2=0$.  Since $\ov R$ is Gorenstein, it follows that
$\ov\fm=0$ or $\mu(\ov\fm)=1$, hence $e=1$ or $e=2$. We may then assume  $\ov \fm^2\ne 0$. If $\fm^3=0$, then, since  $R$ is
Gorenstein, $\fm^2\neq 0$ and $a\in \fm\smallsetminus \fm^2$, it follows that $\fm^2=a\fm$ and thus $\ov\fm^2=0$, a contradiction.

We have thus $\fm^3\ne 0$. Since $R$ is Gorenstein with $\fm^4=0$, we
have $\fm^3=(0:_R\fm)$ and $\mu(\fm^3)=1$. Remark \ref{basic}(4) gives 
\begin{equation}\label{lR}
\lambda(R)=\mu(\fm^2)+e+2.
\end{equation}
The  inclusion
$a\fm^2\subseteq \fm^3$ implies $a\fm^2=\fm^3$ or $a\fm^2=0$. If
$a\fm^2=0$, then $\ds{a\fm\subseteq (0:_R\fm)}$. Since $a\in\fm\smallsetminus
\fm^2$, we have $a\fm\ne 0$, hence $a\fm=(0:_R\fm)=\fm^3$. 

We have thus one of the two cases:
\begin{enumerate}[\quad(i)]
\item  $a\fm=\fm^3$ and $a\fm^2=0$, hence  $\mu(a\fm)=1$ and
$\mu(a\fm^2)=0$. 
\item $a\fm^2=\fm^3$, hence   $\mu(a\fm^2)=1$. 
\end{enumerate}
In either case, $\fm^3\subseteq a\fm$, hence $\ov \fm^3=0\neq\ov\fm^2$. 
 Since $\ov R$ is Gorenstein, we conclude
that $\rank_k(\ov\fm^2)= 1$, hence
$\rank_k(\fm^2/a\fm)=1$ by \eqref{am}. A length count in the exact
sequence
$$
0\to a\fm/\fm^3\to\fm^2/\fm^3\to\fm^2/a\fm\to 0
$$
gives 
\begin{equation}
\label{a}
\rank_k(a\fm/\fm^3)=\mu(\fm^2)-1
\end{equation}

If (i) holds, then $a\fm/\fm^3=0$, hence  $\mu(\fm^2)=1$.  Remark
\ref{basic} gives
\[
\lambda(aR)=\mu(a\fm)+\mu(a\fm^2)+1=2=\mu(\fm^2)+1
\]

If (ii) holds, then $\mu(a\fm)=\mu(\fm^2)-1$ by \eqref{a} and Remark
\ref{basic} gives:
\[
\lambda(aR)=\mu(a\fm)+\mu(a\fm^2)+1
=\mu(\fm^2)+1
\]
In either case, Remark \ref{basic} and \eqref{lR} yield
\[
\lambda(0:a)=\lambda(R)-\lambda(aR)=\big(\mu(\fm^2)+e+2\big)-(\mu(\fm^2)+1)=e+1
\]

Similar formulas hold for $\lambda(bR)$ and $\lambda(0:b)$, hence 
$$
\mu(\fm^2)+1=\lambda(aR)=\lambda(0:b)=e+1
$$
It follows that $\mu(\fm^2)=e$, thus $\HH Rt=1+et+et^2+t^3$.

We now establish the equivalence of conditions {\rm (3.1)} and {\rm (3.2)}  for $a,b$ in  $\fm\smallsetminus \fm^2$. If {\rm (3.1)}  holds, then $\HH R{-1}=0$ by the above argument. The equalities $\mu(a\fm)=e-1=\mu(b\fm)$ then follow from Proposition \ref{1efs:EPZD}. 
Conversely, assume $ab=0$ and $\mu(a\fm)=e-1=\mu(b\fm)$. 
If $a\fm=0$ then $\mu(a\fm)=e-1=0$, contradicting the hypothesis $e\geq 3$. 
If $\fm^3=0$ then, as $a\fm\ne 0$ and $\lambda(0:_{R}\fm)=1$, one has $a\fm=\fm^2=(0:_{R}\fm)$ and $\mu(a\fm)=1$. The equality $\mu(a\fm)=e-1$ gives $e=2$, contradicting the hypothesis. 

We have thus $\fm^3\ne 0$.  Since $R$ is
Gorenstein, one has $\ds{(0:\fm)=\fm^3}$ and $\mu(\fm^3)=1$. Note that
$a\fm^2\subseteq \fm^3$ hence we have  $a\fm^2=0$ or $a\fm^2=\fm^3$. We
will show that that $a\fm^2\ne 0$, implying that $a\fm^2=\fm^3$.
Similarly, $b\fm^2=\fm^3$ and  Proposition \ref{1efs:EPZD} shows that
$a,b$ is an exact pair of zero divisors.

Indeed, assume  $a\fm^2=0$. We then have $a\fm\subseteq (0:\fm)$,
hence $a\fm\subseteq \fm^3$. If  $a\fm=0$, then  $a\in (0:\fm)=\fm^3$, a
contradiction. It remains thus that $a\fm=\fm^3$, and hence
$\mu(a\fm)=1$. The equality $\mu(a\fm)=e-1$ then yields $e=2$, a contradiction.

To establish the equivalence of $(1)$ and $(2)$ it is enough to recall that, since $R$ is Gorenstein, one has $(0:_{R}(0:_{R}I))=I$ for every ideal $I$.

It then remains to establish that $(2)$ implies (3) (i), since the reverse implication is clear. For this, it is enough to show that $a,b$ are in $\fm\smallsetminus\fm^2$. Note that $0\neq a \in\fm$ implies $0\neq b\in\fm$. As $R$ is Gorenstein $\mu(\fm^3)\leq\rank_{k}(0:_{R}\fm)=1$, thus Lemma \ref{EZD:CZD:m2} (2) yields $a,b\in\fm\smallsetminus\fm^2$, concluding the equivalence of (2) and (3).
\end{proof}

\begin{theorem}\label{gorenstein:structural}
Let $(R,\fm,k)$ be a local  Gorenstein ring with $\mu(\fm)\ge 3$ and $\fm^4=0$. 

If there exists a non-zero element $a\in\fm$ such that the ideal $(0:_{R}a)$ is principal, then  the following hold:
\begin{enumerate}[\quad\rm(1)]
\item $\agr R$ is Gorenstein;
\item $R$ is Koszul;
\item every finite $R$-module $M$ has a syzygy module that is Koszul.
\end{enumerate}
\end{theorem}

\begin{proof}Proposition  \ref{Gorenstein:EZD} shows that $\HH{R}{t}=1+et+et^2+t^3$ with $e=\mu(\fm)$ and  yields the existence of a element $b$ such that $a,b\in \fm\smallsetminus\fm^2$ form an exact pair of zero divisors.

From Remark \ref{R/a-Gorenstein}, $\ov R=R/aR$ is Gorenstein with ${\ov{\fm}}^3=0$, thus $(0:_R{\ov\fm})=\ov\fm^2$ and $\mu(\ov\fm^2)=1$. It follows that $(0:_{\agr R}\agr{\ov\fm})=\agr{\ov\fm}^2$ and $\mu(\agr{\ov\fm}^2)=1$, hence $\agr{\ov R}$ is Gorenstein as well. 
 From Theorem \ref{1efs:thm1}(2) and (1), $\agr R/a^*\agr R$ is thus a Gorenstein local ring where $a^*$ is an exact zero divisor. It then follows from  Remark \ref{R/a-Gorenstein} that $\agr R$ is Gorenstein,  establishing (1).

Since $\rank_k(0:_R\fm)=1$,  the hypothesis $\mu(\fm)\ge 3$ establishes condition (2) of Theorem \ref{1ee1:Koszul:Syzygy:Good} concluding (2) and (3).
\end{proof}

\begin{remark}
The set of standard graded Gorenstein $k$-algebras $R$ with $\HH R{t}=1+et+et^2+t^3$ is in bijective correspondence with the set of degree $3$ forms (up to scalars) in $e$ variables over $k$, via the `inverse systems' correspondence of Macaulay, as described in \cite[\S 6]{CRV}. As a consequence, such algebras can be parametrized by means of a projective space $\mathbf P(A_3)$, where $A_3$ is the degree three graded component of a polynomial ring over $k$ in $e$ variables. In view of Remark \ref{connection:CRV}, the proof of  \cite[6.4]{CRV} shows that when $e\geq 3$ there exists a non-empty Zariski open subset of $\mathbf P(A_3)$, whose points correspond to Gorenstein standard graded $k$-algebras admitting an exact pair of zero divisors (cf \cite[6.5]{CRV}).  In other words: a generic Gorenstein standard graded $k$-algebra of socle degree $3$ has an exact pair of zero divisors in degree one. 
\end{remark}

We now derive the numerical version of Theorem \ref{gorenstein:structural} as presented in the introduction. 

\begin{corollary}\label{gorenstein:numerical}
Let $(R,\fm,k)$ be Gorenstein with $\fm^4=0$ and $\mu(\fm)=e\geq3$.

Assume there exist a non-zero element $a\in\fm$ such that the ideal $(0:_{R}a)$ is principal then the following hold:\begin{enumerate}[\quad\rm(1)]
\item $\HH{R}{t}=1+et+et^2+t^3$;
\item $\HH{R}{-t}\Poi k R t=1$;
\item
$\HH{R}{-t}\Poi M R t$ is in $\mathbb{Z}[t]$, for every finite $R$-module $M$.
\end{enumerate}
\end{corollary}

\begin{proof}
(1) It follows directly from Proposition \ref{Gorenstein:EZD} that  $\HH{R}t=1+et+et^2+t^3$.

(2) and (3) follow from the structural conditions (2) and (3) in Theorem \ref{gorenstein:structural} together with \eqref{Poi:Koszul} and  \ref{koszulsyzygymodule} respectively.
\end{proof}

To finish, we discuss the condition on $\mu(\fm)$ in the theorem. Any Gorenstein ring with $\mu(\fm)\le 2$ is a complete intersection, see \cite[Prop. 5]{Serre}, \cite[\S 5]{Bass}, and Poincar\'e series over such rings are well understood. 

\begin{bfchunk}{Artinian complete intersections.}
\label{minimalmult:ci}
Let $(R,\fm,k)$ be an Artinian complete intersection ring, $M$ be a finite module over $R$ and let $\epsilon (R)$ denote the  Hilbert-Samuel multiplicity of $R$. 
One always has: 
\begin{equation}
\label{eq:ci}
\epsilon(R)\ge 2^{\mu(\fm)
}, \ \Poi k R t=(1-t)^{-\mu(\fm)}\ \text{  and  }\  
(1-t^2)^{\mu(\fm)}\Poi M R t\in\mathbb{Z}[t],
\end{equation}
(cf \cite[\S 7, \text{Prop. 7}]{{BourAC8}}, \cite[Thm. 6]{Tate} and \cite[4.1]{GuLe74}).

Furthermore, the following conditions are equivalent:
\begin{enumerate}[\rm \quad(1)]
\item
$R$ has \textit{minimal multiplicity, i.e. $\epsilon(R)= 2^{\mu(\fm)}$},
\item 
$R$ is Koszul,
\item
$\HH R t=(1+t)^{\mu(\fm)}$;
 \end{enumerate}
 and when these conditions hold, one has  
 \begin{equation}
\label{eq:minci}
(1-t)^{\mu(\fm)} \Poi M R t\in \mathbb{Z}[t].
\end{equation}

  Indeed, Avramov shows in  \cite[(2.3)]{Avr94} that $(\rm 1)$ implies $(\rm 2)$ and \eqref{eq:minci}.
  The implication $\rm (2)\Rightarrow\rm (3)$ follows from \eqref{Poi:Koszul} and \eqref{eq:ci}. On the other hand one has $\epsilon(R)=\lambda(R)=\HH R{1}$ 
 as $R$ is Artinian, so it becomes clear that $\rm(3)\Rightarrow\rm(1)$.
\end{bfchunk}

\section*{Acknowledgments}
The authors are grateful to Luchezar Avramov for reading many versions of this manuscript and making suggestions that improved greatly the content and the presentation. They also wish to thank  Aldo Conca for pointing out that generic Gorenstein algebras of socle degree 3 have an exact pair of zero divisors. 



\begin{thebibliography}{99}
	
	\bibitem{Anick}
		D.~J.~Anick, 
		\textit{A counterexample to a conjecture of Serre}, 
		Ann. of Math. {\bf  115} (1982); 1--33.

	\bibitem{Avr94}
		L.~L.~Avramov,
		\textit{Local rings over which all modules have rational Poincar\'e Series}, 
		J. Pure Appl. Algebra {\bf 91} (1994); 29--48.

	\bibitem{AHS}
  	L.~L.~Avramov, I.~B.~Henriques, L.~M.~\c Sega,
		\textit{Exact zero divisors}, 
		in preparation.

	\bibitem{AKM}
		L.~L.~Avramov, A.~R.~Kustin, M.~Miller, 
		\textit{Poincar\'e series of modules over local rings of small embedding codepth or small linking number},  
		J. Algebra  {\bf 118}  (1988); 162--204. 

	\bibitem{AIS}
		L.~L.~Avramov, S.~Iyengar, L.~M.~\c Sega, 
		\textit{Free resolutions over short local rings},
		J. Lond. Math. Soc. {\bf 78} (2008); 459--476.

	\bibitem{Bass}
		H.~Bass,
		\textit{On the ubiquity of Gorenstein rings},
		Math. Zeitschr. {\bf 82} (1963); 8--28.

	\bibitem{Bogvad}
		R.~B\o gvad,
		\textit{Gorenstein rings with transcendental Poincar\'e series},
		Math. Scand. {\bf 53} (1983); 5--15.

	\bibitem{BourAC8}
    N. Bourbaki,
    \textit{Alg\`ebre commutative, Chapitre VIII, Dimension}, 
    Masson, Paris (1983).

	\bibitem{BourAC9}
    N. Bourbaki,
    \textit{Alg\`ebre commutative, Chapitre IX, Anneaux locaux Noeth\'eriens complets}, 
    Masson, Paris (1983).

	\bibitem{CRV}
		A.~Conca, M.~E.~Rossi, G.~Valla, 
		\textit{Gr\"obner Flags and Gorenstein Algebras},
		Compos. Math. {\bf 129} (2001); 95--121.
		
	\bibitem{Fro99}
		R.~Fr\"oberg,
		\textit{Koszul algebras},
		Advances in Commutative Ring Theory, (Fez, 1999), 
		Lect. Notes Pure Appl. Math. {\bf 205}, Marcel Dekker, 
		New York (1999); 337--350.

	\bibitem{GuLe74}
		T.~Gulliksen, 
		\textit{A change of rings theorem with applications to Poincar\'e series and intersection multiplicity}, 
		Math. Scand. {\bf 34} (1974); 167--183. 

	\bibitem{HI}
		J.~Herzog, S.~Iyengar
		\textit{Koszul modules},
		J. Pure Appl. Algebra {\bf 201} (2005); 154--188.

	\bibitem{Jac} 
		C.~Jacobsson, 
		\textit{Finitely presented graded Lie algebras and homomorphisms of local rings},   
		J. Pure Appl. Algebra  {\bf 38}  (1985); 243--253. 

	\bibitem{JKM}
		C.~Jacobsson, A.~R.~Kustin, M.~Miller
		\textit{The Poincar\'e series of a codimension four Gorenstein ring is rational}, 
		J. Pure Appl. Algebra {\bf 38} (1985); 255--275. 

	\bibitem{Lev75}
		G.~Levin,
		\textit{Local rings and Golod homomorphisms},
		J. Algebra {\bf 37} (1975); 266--289.

	\bibitem{Lev80}
		G.~Levin,
		\textit{Large homomorphisms of local rings},
		Math. Scand. {\bf 46} (1980); 209--215.
		
	\bibitem{Nic04}	
	W.~K.~Nicholson,
	\textit{Rings with the dual of the isomorphism theorem}, J. Algebra  {\bf 271}  (2004);  no. 1, 391--406. 
	\bibitem{Roos}
		J.-E.~Roos, 
		\textit{Good and bad Koszul algebras and their Hochschild homology},  
		J. Pure Appl. Algebra  {\bf 201} (2005); 295--327.

	\bibitem{Sega}
		L.~M.~\c Sega,		
		\textit{Homological properties of powers of the maximal ideal}, 
		J. Algebra {\bf 241} (2001); 827--858.

	\bibitem{Serre}
		J.~P.~Serre,
		\textit{Sur les modules projectifs}, 
		S\'eminaire Dubreil, Alg\`ebre et th\'eorie des nombres, 
		\textbf{14} no. 1 (1960-1961), Expos\'e No. 2. 

	\bibitem{SJO79}		
		G.~Sj\"odin,
		\textit{The Poincar\'e series of modules over Gorenstein rings with} $\fm^3=0$, 
		Mathematiska Institutionen, Stockholms Universitet, 
		Preprint {\bf 2} (1979).

	\bibitem{Soto}
		J.~J.~M.~Soto,		
		\textit{Gorenstein quotients by principal ideals of free Koszul homology}, 
		Glasgow Math. J. {\bf 42} (2000); 51--54.


	\bibitem{Tate}
		J.~Tate,
		\textit{Homology of Noetherian rings and of local rings},
		Illinois J. Math. {\bf 1} (1957); 14--27.

	\bibitem{VaVa}
		P.~Valabrega, G.~Valla,
		\textit{Form rings and regular sequences},
		Nagoya Math. J. {\bf 72} (1978); 93--101.
		
		\end{thebibliography}
     \end{document}